\documentstyle[11pt]{article}
\def\Bbb R{{\rm \bf R}}
\def\proclaim#1{\vskip2mm{\bf #1}\em}
\def\endproclaim{\em \vskip2mm}
\def\tag#1{\eqno(#1)}
\def\gathered{\begin{array}{c}}
\def\endgathered{\end{array}}
\def\text{\mbox}

\begin{document}

\title {On reconstruction in the inverse conductivity problem with
one measurement}
\author{Masaru IKEHATA\footnote{
Department of Mathematics,
Faculty of Engineering,
Gunma University, Kiryu 376-8515, JAPAN}}
\date{}
\maketitle

\begin{abstract}
We consider an inverse problem for electrically conductive material occupying a
domain $\Omega$ in $\Bbb R^2$. Let $\gamma$
 be the conductivity of $\Omega$, and $D$ a subdomain of $\Omega$. We assume that
$\gamma$ is a positive constant $k$ on $D$, $k\not=1$ and is $1$  on $\Omega\setminus D$;
both $D$ and $k$ are unknown. The problem
is to find a reconstruction formula of $D$ from the Cauchy data on $\partial\Omega$
 of a non-constant solution $u$
of the equation $\nabla\cdot\gamma\nabla u=0$ in $\Omega$. 
We prove that if $D$ is known to be a convex polygon such that
$\text{diam}\,D<\text{dist}\,(D,\partial\Omega)$, there are two formulae for calculating the support function of 
$D$ from the Cauchy data.

\noindent
AMS: 35R05

$\quad$

\noindent
Key words: Inverse conductivity problem, Exponentially growing solution, Cauchy data,
support function


\end{abstract}


\section{Introduction}

This paper is the sequel to \cite{Ik} and, as predicted therein, we return to one of the problems treated
by  Friedman-Isakov \cite{FI}.  They considered an inverse problem for electrically conductive
material occupying a bounded domain $\Omega$ in $\Bbb R^2$. Let $\gamma$
be the conductivity of $\Omega$, and $D$
a subdomain of $\Omega$ such that  $\overline D\subset\Omega$.
They assume that $\gamma$ is a positive constant $k$ on $D$
with $k\not=1$ and is $1$ on $\Omega\setminus D$.
Let $u$ be a non-constant solution to the equation
$$\displaystyle
\nabla\cdot\gamma\nabla u=0\,\,\text{in}\,\Omega.
\tag {1.1}
$$
Let $\nu$ denote the unit outward normal vector field to $\Omega\setminus D$.

They considered the following uniqueness problem.

{\bf\noindent Uniqueness problem}

\noindent
Assume that $k$ is known and $D$ is unknown. Can one determine $D$ from the Cauchy data 
$u\vert_{\partial\Omega}$, $\frac{\partial u}{\partial\nu}\vert_{\partial\Omega}$?

They proved that if $D$ is known to be a convex polygon such that
$$\displaystyle
\text{diam}\,D<\text{dist}\, (D, \partial\Omega), 
\tag {1.2}
$$
the answer to the problem is yes.

A strong point of their result is that there is no additional assumption on the behaviour of
$u\vert_{\partial\Omega}$ or $\frac{\partial u}{\partial\nu}\vert_{\partial\Omega}$
at the cost of (1.2). Barcelo et al \cite{BFS} added such an assumption and dropped (1.2).
Seo \cite{S} proved a uniqueness theorem from two sets of the Cauchy data having an additional
restriction on the behaviour and removed (1.2) and the convexity restriction on $D$. When
$\partial D$ has a special geometry, there are some results. For example, Kang-Seo \cite{KS}
obtained a uniqueness result when $D$ is a disc.

If both $D$ and $k$ are unknown, the problem becomes more difficult. Alessandrini-Isakov \cite{AI}
considered this problem and obtained a uniqueness theorem of a convex polygon $D$ and $k$
without (1.2). Instead of this assumption they assume that $u\vert_{\partial\Omega}$ or
$\frac{\partial u}{\partial\nu}\vert_{\partial\Omega}$ has a special property.

From these investigations one can say that the Cauchy data of a solution to (1.1) contain
information about the location of $D$. However, their proofs do not tell us how to extract such
information from the Cauchy data.

In this paper we consider the following reconstruction problem.

{\bf\noindent Reconstruction problem}

\noindent
Assume that both $k $ and $D$ are unknown. Find a formula for calculating information about the
location of $D$ from the Cauchy data of $u$.

This is a purely mathematical problem and remains open. In \cite{Ik} we considered the
extreme case $k = 0$, and obtained such formulae provided $D$ was a convex polygon with the
restriction (1.2). In this paper using the idea discovered therein we present such formulae
under the same geometric assumption on $D$ when $k > 0$, $k\not=1$.

Now we describe the result more precisely. Let $S^1$ denote the set of all unit vectors of  $\Bbb R^2$.
Recall the definition of the support function:
$$\displaystyle
h_D(\omega)=\sup_{x\in D}\,x\cdot\omega,\,\,\omega\in S^1.
$$
From this function one can reconstruct the convex hull of general domain $D$.

We say that $\omega\in S^1$ is regular with respect to $D$ if the set
$$\displaystyle
\{x\in\Bbb R^2\,\vert\,x\cdot\omega=h_D(\omega)\}\cap\partial D
$$
consists of only one point.

{\bf\noindent Remark 1.1.}
Note that if $D$ is a polygon, the counting number of the set of all unit vectors
which are not regular with respect to $D$ is finite. Therefore, it is very rare for us to choose a
direction $\omega$ that is not regular with respect to $D$; 
$h_D(\,\cdot\,)$ is a continuous function. Therefore,
the support function of $D$ is uniquely determined by knowing its restriction to the set of all
unit vectors which are regular with respect to $D$.

We merely assume that $\partial\Omega$ is Lipschitz and $u\in H^1(\Omega)$, and consequently we have to
clarify what we mean by the symbol $\frac{\partial u}{\partial\nu}\vert_{\partial\Omega}$.
 It is defined as an element of the dual space of
$H^{1/2}(\partial\Omega)$ by the formula
$$\displaystyle
<\frac{\partial u}{\partial\nu}\vert_{\partial\Omega},\,f>
=\int_{\Omega}
\{1 + (k-1)\chi_D\}\nabla u\cdot\nabla\Psi dx
\tag {1.3}
$$
where $f\in H^{1/2}(\partial\Omega)$, $\Psi$ is in $H^1(\Omega)$
and satisfies $\Psi = f $ on $\partial\Omega$.
From the definition of the
weak solution we know that it is well defined and one may take $\Psi$ such that $\Psi(x)=0$ 
for $x$ far from  $\partial\Omega$. This means that 
$\frac{\partial u}{\partial\nu}\vert_{\partial\Omega}$ is uniquely determined by the value of $u$ near
$\partial\Omega$.   
We call $ (u\vert_{\partial\Omega}, \frac{\partial u}{\partial\nu}\vert_{\partial\Omega})$
the Cauchy data of $u$ on $\partial\Omega$. It is a pair of the voltage potential and electric
current distribution on $\partial\Omega$.

In this paper the following special harmonic functions are extremely important:
$$\displaystyle
v=v(x)=e^{\tau x\cdot(\omega+i\omega^{\perp})},\,\,\tau>0
$$
where $\omega$, $\omega^{\perp}\in S^1$ and satisfy
$$\displaystyle
\omega\cdot\omega^{\perp}=0,\,\,\text{det}\,(\omega\,\,\omega^{\perp})<0.
$$

{\bf\noindent Remark 1.2.} 
Calder\'on \cite{C} made use of these types of harmonic functions in the inverse
conductivity problem with infinitely many measurements.

Using these functions and the Cauchy data of $u$ on $\partial\Omega$ we give the following definition.

{\bf\noindent Definition 1.1 (Indicator function).}
Let $u$ be a weak solution to (1.1). Define
$$\displaystyle
I_{\omega}(\tau,t)=e^{-\tau t}
\left\{<\frac{\partial u}{\partial\nu}\vert_{\partial\Omega}, v\vert_{\partial\Omega}>
-<\frac{\partial v}{\partial\nu}\vert_{\partial\Omega}, u\vert_{\partial\Omega}>\right\},
\,\,\tau>0,\,\,t\in\Bbb R.
$$
Note that $u$ is fixed. The result is the two following formulae.

\proclaim{\noindent Theorem 1.1.}
Assume that $D$ is a convex polygon satisfying  (1.2) and that $u$ is not a constant
function. Let $\omega$ be regular with respect to $D$. The formulae
$$\displaystyle
\{t\in\Bbb R\,\vert\,\lim_{\tau\longrightarrow\infty}I_{\omega}(\tau,t)=0\}
=[h_D(\omega),\,\infty[
\tag {1.4}
$$
$$
h_D(\omega)-t
=\lim_{\tau\longrightarrow\infty}
\frac{\log\vert I_{\omega}(\tau,t)\vert}{\tau},\,\,\forall t\in\Bbb R,
\tag {1.5}
$$
are valid.

\endproclaim

This is a direct corollary of the trivial identity
$$\displaystyle
I_{\omega}(\tau,t)=e^{\tau(h_D(\omega)-t)}I_{\omega}(\tau, h_D(\omega))
$$
and the asymptotic behaviour of $I_{\omega}(\tau, h_D(\omega))$ as $\tau\longrightarrow\infty$
described below.

\proclaim{\noindent Key lemma.}
Assume that $D$ is a convex polygon satisfying (1.2) and that $u$ is not a constant
function. Let $\omega$ be regular with respect to $D$. There exist positive constants $L$ and $\mu$ such that
$$\displaystyle
\lim_{\tau\longrightarrow\infty}\tau^{\mu}\vert I_{\omega}(\tau,h_D(\omega))\vert=L.
$$
\endproclaim

The proof of this lemma is delicate and the outline is as follows. From the regularity of
$\omega$ we know that the line $x\cdot\omega=h_D(\omega)$ meets 
$\partial D$ at a vertex $x_0$ of $D$.
Using a well known
expansion of $u$ about $x_0$ (see proposition 2.1) and a formula which connects 
$I_{\omega}(\tau, h_D(\omega))$ with
an integral on $\partial D$ involving $u\vert_{\partial D}$ (see proposition 3.1), 
we obtain the asymptotic expansion of
$I_{\omega}(\tau,h_D(\omega))$ as $\tau\longrightarrow\infty$ (see proposition 3.2):
$$\displaystyle
I_{\omega}(\tau,h_D(\omega))\sim e^{i\tau x_0\cdot\omega}
\sum_{j=1}^{\infty}\frac{L_j}{\tau^{\mu_j}}
$$
where $0<\mu_1<\mu_2<\cdots$.  The problem is to show that $L_j\not=0$ for some $j$ . 
We see that if
$L_j=0$ for all $j$ , $u$ has a harmonic continuation in a neighbourhood of $x_0$ (see lemma 4.1).
Then Friedman-Isakov's  extension argument \cite{FI} tells us that $u$ has to be a constant function
and it is a contradiction. Restriction (1.2) is merely employed to make use of their argument.

It would be interesting to apply our method to the three-dimensional problem (see \cite{BFS}
for a uniqueness result) or a similar problem in the linear theory of elasticity. This will be
considered in subsequent papers. The numerical testing of (1.4) and (1.5) remains open and
we hope that someone performs this task in the future.

Finally, we note that in subsequent sections we always assume that $\omega$ is regular with
respect to $D$.

\section{Preliminaries}

\subsection{Notation}

\noindent
$x_0$ stands for the only one point of the set
$$\begin{array}{c}
\displaystyle
\{x\in\Bbb R^2\,\vert\,x\cdot\omega=h_D(\omega)\}\cap\partial D;
\\
\\
\displaystyle
B_R(x_0)=\{x\in\Bbb R^2\,\vert\,\vert x-x_0\vert<R\},\,\,R>0;
\end{array}
$$
$\Theta$ stands for the outside angle at the vertex $x_0$ of $D$ and thus $\pi<\Theta<2\pi$.

\subsection{Expansion of u about a vertex}

\noindent
Let $u$ be a weak solution to (1.1). Define
$$\begin{array}{c}
\displaystyle
u^e=u\vert_{\Omega\setminus D}\\
\\
\displaystyle
u^i=u\vert_D.
\end{array}
$$
We introduce polar coordinates. Let $\omega^{\perp}$ denote the unit-vector perpendicular to $\omega$
satisfying $\text{det}\,(\omega\,\,\omega^{\perp})<0$.
Since $x_0$ is vertex of $D$ and $\omega$ is regular with respect to $D$, one may write
$$\begin{array}{c}
\displaystyle
B_{2\eta}(x_0)\cap (\Omega\setminus\overline D)
=\{x_0+r(\cos\theta\,a+\sin\theta\,a^{\perp})\,\vert\,
0<r<2\eta,\,0<\theta<\Theta\}\\
\\
\displaystyle
B_{2\eta}(x_0)\cap\overline D=\{x_0+r(\cos\theta\,a+\sin\theta\,a^{\perp})\,\vert\,
0<r<2\eta,\,\Theta<\theta<2\pi\}
\\
\\
\displaystyle
B(x_0,\eta)\cap\partial D=\Gamma_p\cup\Gamma_q\cup\{x_0\}\\
\\
\displaystyle
\Gamma_p=\{x_0+r(\cos p\,\omega^{\perp}+\sin p\,\omega)\,\vert\,0<r<\eta\}\\
\\
\displaystyle
\Gamma_q=\{x_0+r(\cos q\,\omega^{\perp}+\sin q\,\omega)\,\vert\,0<r<\eta\}
\end{array}
$$
where $\eta$ is a small positive number,
$$\begin{array}{c}
\displaystyle
-\pi<q<p<0\\
\\
\displaystyle
p+\Theta=2\pi+q\\
\\
\displaystyle
a=\cos p\,\omega^{\perp}+\sin p\,\omega\\
\\
\displaystyle
a^{\perp}=-\sin p\,\omega^{\perp}+\cos p\,\omega\\
\\
\displaystyle
\text{det}\,(a\,\,a^{\perp})>0.
\end{array}
$$
Set
$$\displaystyle
u(r,\theta)=u(x),\,\,x=x_0+r(\cos \theta\,a+\sin \theta\,a^{\perp}).
$$
The following proposition is only given for our purpose and the proof is well known. For
example, the reader can find its outline in section 2 of \cite{BFI}.

\proclaim{\noindent Proposition 2.1.}
There exist a real number $\alpha$, a monotone increasing sequence $(\mu_j)_{j=1,\cdots}$ of
positive numbers and sequences 
$\{A^e_j\}$, $\{B^e_j\}$, $\{A^i_j\}$, $\{B^i_j\}$ of real numbers such that:
$$\displaystyle
(1+k)^2\sin^2\pi\mu_j=
(1-k)^2\sin^2(\pi-\Theta)\mu_j;
\tag {2.1}
$$
$$\displaystyle
\left(\begin{array}{c}
\displaystyle
A^e_j\\
\\
\displaystyle
B^e_j
\end{array}
\right)
=\left(\begin{array}{lr}
\displaystyle
\cos 2\pi\mu_j & \displaystyle \sin 2\pi\mu_j\\
\\
\displaystyle
-k\sin 2\pi\mu_j & \displaystyle k\cos 2\pi\mu_j
\end{array}
\right)
\left(\begin{array}{c}
\displaystyle
A^i_j\\
\\
\displaystyle
B^i_j
\end{array}
\right),
\tag {2.2}
$$
$$\displaystyle
\left(\begin{array}{c}
\displaystyle
A^e_j\\
\\
\displaystyle
B^e_j
\end{array}
\right)
=\left(\begin{array}{lr}
\displaystyle
\cos^2\Theta\mu_j+k\sin^2\Theta\mu_j & \displaystyle 
(1-k)\cos\Theta\mu_j\sin\Theta\mu_j\\
\\
\displaystyle
(1-k)\cos\Theta\mu_j\sin\Theta\mu_j & \displaystyle 
\sin^2\Theta\mu_j+k\cos^2\Theta\mu_j
\end{array}
\right)
\left(\begin{array}{c}
\displaystyle
A^i_j\\
\\
\displaystyle
B^i_j
\end{array}
\right),
\tag {2.3}
$$
$$\begin{array}{c}
\displaystyle
u^e(r,\theta)-\alpha=\sum_{j=1}^{\infty}r^{\mu_j}(A^e_j\cos\mu_j\theta+B^e_j\sin\mu_j\theta),\\
\\
\displaystyle
u^i(r,\theta)-\alpha=\sum_{j=1}^{\infty}r^{\mu_j}(A^i_j\cos\mu_j\theta+B^i_j\sin\mu_j\theta);
\end{array}
\tag {2.4}
$$
the series are absolutely convergent in $H^1(B_{s\eta}(x_0)\cap(\Omega\setminus\overline D))$
and $H^1(B_{s\eta}(x_0)\cap D)$, respectively, and uniformly in 
$B_{s\eta}(x_0)$ for each $0<s<2$; moreover for each $l=1,\cdots$,
$$\begin{array}{c}
\displaystyle
\left\vert u(r,0)-\alpha-\sum_{j=1}^l r^{\mu_j}A^e_j\right\vert\le C_lr^{\mu_{l+1}}\\
\\
\displaystyle
\left\vert u(r,\Theta)-\alpha-\sum_{j=1}^l
r^{\mu_j}(A^i_j\cos\Theta\mu_j+B^i_j\sin\Theta\mu_j)\right\vert\le C_lr^{\mu_{l+1}},
\,\, 0<r<\eta.
\end{array}
\tag {2.5}
$$

\endproclaim

Note that from (2.2) and (2.3) we have

$$\begin{array}{c}
\displaystyle
A^i_j(\cos 2\pi\mu_j-\cos^2\Theta\mu_j-k\sin^2\Theta\mu_j)\\
\\
\displaystyle
+B^i_j\{\sin 2\pi\mu_j+(k-1)\cos\Theta\mu_j\sin\Theta\mu_j\}=0.
\end{array}
\tag {2.6}
$$

\section{Asymptotic expansion of the indicator function}

\proclaim{\noindent Proposition 3.1.}
Let $v$ be a $H^2(\Omega)$ harmonic function.  For any constant $\lambda$ the formula
$$\displaystyle
<\frac{\partial u}{\partial\nu}\vert_{\partial\Omega}, v\vert_{\partial\Omega}>
-<\frac{\partial v}{\partial\nu}\vert_{\partial\Omega},u\vert_{\partial\Omega}>
=
(1-k)\int_{\partial D}(u-\lambda)\frac{\partial v}{\partial\nu},
\tag {3.1}
$$
is valid.

\endproclaim

{\it\noindent Proof.}
From (1.3) we have
$$\begin{array}{c}
\displaystyle
<\frac{\partial u}{\partial\nu}\vert_{\partial\Omega}, v\vert_{\partial\Omega}>
=\int_{\Omega}\{1+(k-1)\chi_D\}\nabla u\cdot\nabla v dx\\
\\
\displaystyle
<\frac{\partial v}{\partial\nu}\vert_{\partial\Omega}, u\vert_{\partial\Omega}>
=\int_{\Omega}\nabla v\cdot\nabla u dx.
\end{array}
\tag {3.2}
$$
Green's formula (see \cite{G})yields
$$\displaystyle
\int_D\nabla u\cdot\nabla v dx=-\int_{\partial D}(u-\lambda)\frac{\partial v}{\partial\nu}.
\tag {3.3}
$$
Note that $\nu$ is outward to $\Omega\setminus\overline D$.
A combination of (3.2) and (3.3) gives (3.1).

\noindent
$\Box$

\proclaim{\noindent Proposition 3.2.}
The asymptotic expansion
$$\displaystyle
I_{\omega}(\tau, h_D(\omega))
\sim
(k-1)ie^{i\tau x_0\cdot\omega^{\perp}}
\sum_{j=1}^{\infty}
e^{i\frac{\pi}{2}\mu_j}\Gamma(1+\mu_j)K_j\tau^{-\mu_j},
\tag {3.4}
$$
is valid where
$$\displaystyle
K_j=A^e_je^{ip\mu_j}
-(A^i_j\cos\Theta\mu_j+B^i_j\sin\Theta\mu_j)e^{iq\mu_j}.
$$

\endproclaim

{\it\noindent Proof.}
For $\eta$ in Section 2 take a positive constant $c$ in such a way that
$$\displaystyle
\partial D\setminus B_{\eta}(x_0)
\subset
\{x\cdot\omega\le h_D(\omega)-c\}.
$$
It follows from (3.1) that
$$\begin{array}{c}
\displaystyle
\frac{I_{\omega}(\tau,h_D(\omega))}
{1-k}
=e^{-\tau h_D(\omega)}\int_{\partial D}(u-\alpha)\frac{\partial v}{\partial\nu}\\
\\
\displaystyle
=e^{-\tau h_D(\omega)}\int_{\Gamma_p}(u-\alpha)\frac{\partial v}{\partial\nu}
+e^{-\tau h_D(\omega)}\int_{\Gamma_q}(u-\alpha)\frac{\partial v}{\partial\nu}
+O(\tau e^{-c\tau}).
\end{array}
\tag {3.5}
$$
Since
$$\begin{array}{c}
\displaystyle
\nu=\sin p\,\omega^{\perp}-\cos p\,\omega\,\,\text{on}\,\Gamma_p\\
\\
\displaystyle
\nu=-\sin q\,\omega^{\perp}+\cos q\,\omega\,\,\text{on}\,\Gamma_q\\
\\
\displaystyle
x\cdot\omega=h_D(\omega)+r\sin(\theta+p)\\
\\
\displaystyle
x\cdot\omega^{\perp}=x_0\cdot\omega^{\perp}+r\cos(\theta+p)\\
\\
\displaystyle
\nabla v=\tau(\omega+i\omega^{\perp})e^{\tau(x\cdot\omega+ix\cdot\omega^{\perp})},
\end{array}
$$
we have
$$\begin{array}{c}
\displaystyle
e^{-\tau h_D(\omega)}\frac{\partial v}{\partial\nu}
=-\tau e^{-ip}e^{i\tau x_0\cdot\omega^{\perp}}e^{r\tau(\sin p+i\cos p)}\,\,\text{on}\,\Gamma_p\\
\\
\displaystyle
e^{-\tau h_D(\omega)}\frac{\partial v}{\partial\nu}
=\tau e^{-iq}e^{i\tau x_0\cdot\omega^{\perp}}e^{r\tau(\sin q+i\cos q)}\,\,\text{on}\,\Gamma_q.
\end{array}
\tag {3.6}
$$
From (2.5) and (3.6) we obtain
$$\begin{array}{c}
\displaystyle
e^{-\tau h_D(\omega)}
\int_{\Gamma_p}\left(u-\alpha-\sum_{j=1}^l r^{\mu_j}A^e_j\right)\frac{\partial v}{\partial\nu}
=O\left(\frac{1}{\tau^{\mu_{l+1}}}\right),
\\
\\
\displaystyle
e^{-\tau h_D(\omega)}
\int_{\Gamma_q}\left\{u-\alpha-\sum_{j=1}^l r^{\mu_j}
(A^i_j\cos\Theta\mu_j+B^i_j\sin\Theta\mu_j)\right\}\frac{\partial v}{\partial\nu}
=O\left(\frac{1}{\tau^{\mu_{l+1}}}\right).
\end{array}
\tag {3.7}
$$
A combination of (3.5)-(3.7) gives
$$\begin{array}{c}
\displaystyle
\frac{I_{\omega}(\tau,h_D(\omega))}
{1-k}
=-\tau e^{-ip}e^{i\tau x_0\cdot\omega^{\perp}}
\sum_{j=1}^l A^e_j\int_0^{\eta}r^{\mu_j}e^{r\tau(\sin p+i\cos p)}dr\\
\\
\displaystyle
+\tau e^{-iq}e^{i\tau x_0\cdot\omega^{\perp}}
\sum_{j=1}^l(A^i_j\cos\Theta\mu_j+B^i_j\sin\Theta\mu_j)
\int_0^{\eta}r^{\mu_j}e^{r\tau(\sin q+i\cos q)}dr\\
\\
\displaystyle
+O\left(\frac{1}{\tau^{\mu_{l+1}}}\right).
\end{array}
\tag {3.8}
$$
We make use of the following formulae \cite{Ik}:
$$\begin{array}{c}
\displaystyle
\int_0^{\eta}r^{\mu_j}e^{r\tau(\sin p+i\cos p)}dr
=\tau^{-(1+\mu_j)}
ie^{i\frac{\pi}{2}\mu_j}e^{ip}e^{ip\mu_j}\Gamma(1+\mu_j)
+O\left(\frac{e^{\eta\tau\sin p}}{\tau}\right),\\
\\
\displaystyle
\int_0^{\eta}r^{\mu_j}e^{r\tau(\sin q+i\cos q)}dr
=\tau^{-(1+\mu_j)}
ie^{i\frac{\pi}{2}\mu_j}e^{iq}e^{iq\mu_j}\Gamma(1+\mu_j)
+O\left(\frac{e^{\eta\tau\sin q}}{\tau}\right).
\end{array}
\tag {3.9}
$$
From (3.8) and (3.9) we obtain (3.4).

\noindent
$\Box$

\section{Proof of the key lemma}

\noindent
The problem is: what happens when
$$\displaystyle
K_j=A^e_je^{ip\mu_j}
-(A^i_j\cos\Theta\mu_j+B^i_j\sin\Theta\mu_j)e^{iq\mu_j}=0
$$
for all $j=1,\cdots$ ?

Sine 
$$\displaystyle
p+\Theta=2\pi+q,
$$
we have
$$\displaystyle
e^{i\Theta\mu_j}e^{ip\mu_j}
=e^{i2\pi\mu_j}e^{iq\mu_j}.
$$
So $K_j=0$ if and only if
$$\displaystyle
A^e_je^{i2\pi\mu_j}
=(A^i_j\cos\Theta\mu_j+B^i_j\sin\Theta\mu_j)e^{i\Theta\mu_j}.
\tag {4.1}
$$
Since $A^e_j$, $B^e_j$, $A^i_j$, $B^i_j$ are all real, we know that (4.1) is equivalent to
$$\displaystyle
A^i_j\cos\Theta\mu_j\,\cos(\Theta-2\pi)\mu_j
+B^i_j\sin\Theta\mu_j\,\cos(\Theta-2\pi)\mu_j
=A^e_j
\tag {4.2}
$$
and
$$\displaystyle
A^i_j\cos\Theta\mu_j\,\sin(\Theta-2\pi)\mu_j
+B^i_j\sin\Theta\mu_j\,\sin(\Theta-2\pi)\mu_j=0.
\tag {4.3}
$$
In this section we only consider $j$ satisfying
$$\displaystyle
\left(\begin{array}{c}
\displaystyle
A^i_j\\
\\
\displaystyle
B^i_j
\end{array}
\right)\not=
\left(
\begin{array}{c}
\displaystyle
0\\
\\
\displaystyle
0
\end{array}
\right).
$$
Since $A^i_j$ and $B^i_j$ are non-trivial solutions of (2.6) and (4.3), we obtain
$$\begin{array}{c}
\displaystyle
L\equiv
(\cos 2\pi\mu_j-\cos^2\Theta\mu_j-k\sin^2\Theta\mu_j)
\sin\Theta\mu_j\,\sin(\Theta-2\pi)\mu_j\\
\\
\displaystyle
-\{\sin 2\pi\mu_j+(k-1)\cos\Theta\mu_j\,\sin\Theta\mu_j\}
\cos\Theta\mu_j\,\sin(\Theta-2\pi)\mu_j=0.
\end{array}
\tag {4.4}
$$
Since
$$\begin{array}{c}
\displaystyle
L=\sin(\Theta-2\pi)\mu_j\times\\
\\
\displaystyle
\{\cos 2\pi\mu_j\,\sin\Theta\mu_j-\cos^2\Theta\mu_j\,\sin\Theta\mu_j-k\sin^3\Theta\mu_j
-\sin2\pi\mu_j\,\cos\Theta\mu_j-(k-1)\cos^2\Theta\mu_j\,\sin\Theta\mu_j\}\\
\\
\displaystyle
=\sin(\Theta-2\pi)\mu_j\,\{\sin(\Theta-2\pi)\mu_j-k\sin\Theta\mu_j\}\\
\\
\displaystyle
=\sin(2\pi-\Theta)\mu_j\,\{\sin(2\pi-\Theta)\mu_j+k\sin\Theta\mu_j\}.
\end{array}
$$
Therefore, (4.4) becomes
$$\displaystyle
\sin(2\pi-\Theta)\mu_j\,\{\sin(2\pi-\Theta)\mu_j+k\sin\Theta\mu_j\}=0.
\tag {4.5}
$$
Moreover, from (4.2) and (4.3) it is easy to see that
$$\displaystyle
A^e_j\sin\Theta\mu_j\,\sin(2\pi-\Theta)\mu_j=0.
\tag {4.6}
$$
This is a compatibility condition of the system (4.2) and (4.3).  Now we are ready to prove the central part of this paper.

\proclaim{\noindent Lemma 4.1.}  Assume that $K_j=0$ for all $j=1,\cdots$.
There exist an integer $a\ge 2$ independent of $j$ and a harmonic continuation $\tilde{u}$ of $u$ 
from $\Omega\setminus\overline D$ into $(\Omega\setminus\overline D)\cup B_{\eta}(x_0)$ such that
$$\displaystyle
\tilde{u}\left(r,\theta+\frac{2\pi}{a}\right)=\tilde{u}(r,\theta)\,\,\text{in}\,B_{\eta}(x_0).
$$

\endproclaim

{\it\noindent Proof.}
The proof is divided into three parts.

$\quad$

{\it\noindent Step 1:}  $\sin(2\pi-\Theta)\mu_j=0$.

To prove this we assume that $\sin(2\pi-\Theta)\mu_j\not=0$.  From (4.5) we get
$$\displaystyle
\sin(2\pi-\Theta)\mu_j+k\sin\Theta\mu_j=0
\tag {4.7}
$$
and this thus yields $\sin\Theta\mu_j\not=0$.  From (4.6) we conclude that $A^e_j=0$.  
Then taking the first components of (2.2) and (2.3), respectively, we get
$$\displaystyle
\left(\begin{array}{lr}
\displaystyle \cos 2\pi\mu_j   & \displaystyle \sin 2\pi\mu_j\\
\\
\displaystyle
\cos^2\Theta\mu_j+k\sin^2\Theta\mu_j & \displaystyle (1-k)\cos\Theta\mu_j\,\sin\Theta\mu_j
\end{array}
\right)
\left(\begin{array}{c}
\displaystyle
A^i_j\\
\\
\displaystyle
B^i_j
\end{array}
\right)
=
\left(\begin{array}{c}
\displaystyle
0\\
\\
\displaystyle
0
\end{array}
\right).
$$
Since $A^i_j$, $B^i_j$ are not trivial solutions of this system, we obtain
$$\begin{array}{c}
\displaystyle
0=\cos 2\pi\mu_j\,(1-k)\,\cos\Theta\mu_j\,\sin\Theta\mu_j
-(\cos^2\Theta\mu_j+k\sin^2\Theta\mu_j)\sin2\pi\mu_j\\
\\
\displaystyle
=\cos\Theta\mu_j\,(\cos 2\pi\mu_j\,\sin\Theta\mu_j-\cos\Theta\mu_j\,\sin 2\pi\mu_j)\\
\\
\displaystyle
-k\sin\Theta\mu_j\,(\cos 2\pi\mu_j\,\cos\Theta\mu_j+\sin\Theta\mu_j\,\sin 2\pi\mu_j)\\
\\
\displaystyle
=-(\cos\Theta\mu_j\,\sin(2\pi-\Theta)\mu_j+k\sin\Theta\mu_j\,\cos(2\pi-\Theta)\mu_j).
\end{array}
\tag {4.8}
$$
A combination of 84.7) and (4.8) gives
$$\displaystyle
\sin(2\pi-\Theta)\mu_j\,\{\cos\Theta\mu_j-\cos(2\pi-\Theta)\mu_j\}=0
$$
and this thus yields
$$\displaystyle
\cos\Theta\mu_j=\cos(2\pi-\Theta)\mu_j.
$$
Therefore, we obtain
$$\displaystyle
\vert\sin\Theta\mu_j\vert=\vert\sin(2\pi-\Theta)\mu_j\vert.
\tag {4.9}
$$
A combination of (4.7) and (4.9) yields
$$\begin{array}{c}
\displaystyle
\vert\sin\Theta\mu_j\vert
=\vert\sin(2\pi-\Theta)\mu_j\vert\\
\\
\displaystyle
=k\vert\sin\Theta\mu_j\vert
\end{array}
$$
and hence $k=1$.  This is a contradiction.

$\quad$

{\it\noindent Step 2:}  $\mu_j$ has to be an integer.

It follows from Step 1 that $(2\pi-\Theta)\mu_j=n\pi$ for an integer $n$.
Then $(\pi-\Theta)\mu_j=-\pi\mu_j+n\pi$.  This gives
$$\displaystyle
\sin(\pi-\Theta)\mu_j=(-1)^{n+1}\sin\pi\mu_j.
$$
Combining this with (2.1), we obtain
$$\begin{array}{c}
\displaystyle
(1+k)^2\sin^2\pi\mu_j=(1-k)^2\sin^2(\pi-\Theta)\mu_j\\
\\
\displaystyle
=(1-k)^2\sin^2\pi\mu_j.
\end{array}
$$
Since $k\not=0$, we have the desired conclusion.

$\quad$

{\it\noindent Step 3:}  From Step 1 we know that there exits an integer $n_j$ such that $(2\pi-\Theta)\mu_j=n_j\pi$.
Since $\mu_j\not=0$, we have
$$\displaystyle
\frac{\Theta}{\pi}=2-\frac{n_j}{\mu_j}.
\tag {4.10}
$$
From Step 2 one it concludes that $\frac{\Theta}{\pi}$ has to be a rational number.  Since $\pi<\Theta<2\pi$, one may write
$$\displaystyle
\frac{\Theta}{\pi}=1+\frac{b}{a}
$$
where $a=2,\cdots$, $b=1,\cdots$ with $(a,b)=1$.  Note that $a$ and $b$ are independent of $j$.
From (4.10) and (4.11) we get
$$\displaystyle
b\mu_j=a(\mu_j-n_j).
$$
Since $(a,b)=1$, there exists an integer $l_j$ such that
$$\displaystyle
\mu_j=l_ja.
$$
Then
$$\displaystyle
\left(\theta+\frac{2\pi}{a}\right)\mu_j=\theta\mu_j+2l_j\pi,
\tag {4.12}
$$
and we have
$$\displaystyle
u(r,\theta)
=\alpha+\sum_{j=1}^{\infty}r^{\mu_j}(A^e_j\cos\theta\mu_j+B^e_j\sin\theta\mu_j)\,\,\text{in}\,(\Omega\setminus\overline D)\cap B_{\eta}(x_0).
$$
By virtue of (4.12), this right-hand side gives a desired harmonic continuation of $u$.

\noindent
$\Box$

Now we are ready to prove the key lemma.  Assume that $K_j=0$ for all $j=1,\cdots$.
From a combination of Lemma 4.1 and Friedman-Isakov's extension argument (see proof of Theorem 1.1 on p.570 in \cite{FI})
we obtain a harmonic extension of $u$ into whole $\Omega$.  This yields that $u$ has to be constant.  This is a contradiction.

So one can take 
$$\displaystyle
m=\min\,\{j\,\vert\,K_j\not=0\}.
$$
Then from (3.4) we have
$$\displaystyle
I_{\omega}(\tau,h_D(\omega))
\sim (k-1)ie^{i\tau x_0\cdot\omega^{\perp}}
e^{i\frac{\pi}{2}\mu_m}\Gamma(1+\mu_m)K_m\tau^{-\mu_m}.
$$
This completes the proof.

$$\quad$$

\centerline{{\bf Acknowledgment}}

The author thanks the referees for several suggestions for the improvement of the manuscript.
This research was partially supported by Grant-in-Aid for Scientific Research (Grant
no 11640151), Ministry of Education, Science and Culture, Japan.

$$\quad$$

\end{document}